\documentclass[10pt]{amsart}

\usepackage[latin1]{inputenc}
\usepackage{amsmath}
\usepackage{amsfonts}
\usepackage{amssymb}
\usepackage{ytableau}
\usepackage[mathscr]{eucal}
\usepackage{diagrams}
\usepackage{xy,color,graphicx}
\usepackage{hyperref}
\xyoption{all}
\usepackage{array}
\usepackage{enumerate}
\usepackage{url}
\usepackage{halloweenmath}

\newcommand{\N}{\mathbb{N}}
\newcommand{\Q}{\mathbb{Q}}
\newcommand{\C}{\mathbb{C}}
\newcommand{\Z}{\mathbb{Z}}

\newcommand{\Oscr}{\mathcal{O}}

\newcommand{\F}{\mathbb{F}}

\title{Counting in times of fake fields}
\author{Lieven Le Bruyn} 
\address{Department Mathematics, University  of Antwerp ,
 Middelheimlaan 1, B-2020 Antwerp (Belgium) {\tt lieven.lebruyn@uantwerpen.be}}

\begin{document}
\sloppy

\maketitle

\begin{abstract}
These are notes of a talk, given at the {\em arithm\'etique en plat pays} meeting in february 2020, on the potential uses of geometries over the  fake field $\F_1$ to zeta functions and counting measures on motives.
\end{abstract}

\section{Counting}

The counting problems we will consider here are all of a geometric nature. Let us start with some examples:
\begin{enumerate}
\item{Let $M$ be a manifold, equipped with a diffeomorphism $f$. Then we can consider $(M,f)$ as a {\em discrete dynamical system}: at time $n \in \N$ we consider the diffeomeorphism $f^n$ sending the point $x \in M$ to $f^n(x)=f \circ \hdots \circ f(x)$.
In such situations, it is important to study the {\em periodic orbits}, or equivalently, to count the number of fixed points $\# Fix(M,f^n)$ of $f^n$. If we are in a situation that all these numbers are finite, we can package these numbers together in a {\em zeta-function}, the so called {\em Artin-Mazur zeta function}
\[
\zeta_{AM}(M,f) = exp(\sum_{n=1}^{\infty} \frac{\# Fix(M,f^n)}{n} t^n) \]
and investigate its properties, for example, when this zeta function is a rational function.}
\item{Let $X$ be an algebraic variety defined over the finite field $\F_p$, that is, locally the points of $X$ are the solutions in the algebraic closure $\overline{\F_p}$ of a set of polynomial equations with coefficients in $\F_p$. Here, we are interested in the number of points $\# X(\F_{p^n})$ of $X$ having all there coordinates lying in the finite field $\F_{p^n}$. One again, we can package these (necessarily) finte numbers into the {\em Weil zeta function} of $X$
\[
\zeta_W(X) = exp(\sum_{n=1}^{\infty} \frac{\# X(\F_{p^n})}{n} t^n) \]
There's more than a superficial relation with the manifold case above. If $Fr$ is the {\em Frobenius automorphism} on $X$, which raises all coefficients of a point $x \in X$ to the $p$-th power, then we have
\[
X(\F_{p^n}) = Fix(X(\overline{\F_p}),Fr^n) \]
}
\end{enumerate}

In both cases these numbers satisfy additive and multiplicative properties. For example, in the $(M,f)$ case, if $N$ is a submanifold stable under $f$,then clearly
\[
\# Fix(M,f^n) = \# Fix(M-N,f^n)+ \# Fix(N,f^n) \]
and if $(M',f')$ is another discrete dynamical system, then so is their  product system  $(M \times M',(f,f'))$ and we have that
\[
\# Fix(M \times M',(f,f')^n) = \# Fix(M,f^n) \times \# Fix(M',f'^n) \]
Similarly, if $Y$ is an algebraic subvariety of $X$ defined over $\F_p$, then clearly
\[
\# X(\F_{p^n}) = \# (X-Y)(\F_{p^n}) + \# Y(\F_{^n}) \]
and for another algebraic variety $X'$ defined over $\F_{p^n}$, we have for the product
\[
\# (X \times Y)(\F_{p^n}) = \# X(\F_{p^n}) \times \# X(\F_{p^n}) \]
Also, if $X$ and $X'$ are isomorphic $\F_p$-varieties, then $\# X(\F_{p^n}) = \# X'(\F_{p^n})$, and if $(M'f)$ and $(M',f')$ are isomorphic discrete dynamical systems (there is a diffeomorphism $\phi : M \rTo M'$ such that $f' = \phi \circ f \circ \phi^{-1}$), then $\# Fix(M,f^n) = \# Fix(M',f'^n)$.

\vskip 3mm

We can now formalise what me mean by a {\em geometric counting problem}. Assume we have a class of geometric objects and a suitable notion of isomorphism, behaving well under sub-objects and products. That is, we can define a commutative ring of {\em motives} $\mathbf{Mot}$ whose elements are the isomorphism classes $[ X ]$ of our objects with addition and multiplication defined by
\[
[X] = [X-Y]+[Y] \quad \text{and} \quad [X \times X'] = [X] \times [X'] \]
whenever $Y$ is a subobject of $X$. 

A {\em counting measure} $\mu$ is then the assignment of an integer $\mu(X) \in \Z$ to any object $X$, invariant under isomorphism, such that the induced map
\[
\mu~:~\mathbf{Mot} \rTo \Z \]
is a ringmorphism.

\section{Fake fields}

Our geometric counting is done in the ring of integers $\Z$ which is a very special object in $\mathbf{comm}$, the category of all commutative rings, in that it is the {\em initial object} of $\mathbf{comm}$, that is, for any commutative ring $R$, there is a unique ringmorphism $i_R~:~\Z \rTo R$.

The geometric object corresponding to the commutative ring $R$ is its {\em affine scheme} $\mathbf{Spec}(R)$, which is the topological space of all prime ideals of $R$ (with the Zariski topology) together with the structure sheaf $\Oscr_R$ on it. As the direction of arrows reverses ingoing from commutative rings to affine schemes this means that the geometric object corresponding to $\Z$, that is $\mathbf{Spec}(\Z)$, is the {\em terminal object} in schemes, that is, any scheme $X$ maps to it and as such one would expect $\mathbf{Spec}(\Z)$ to be the simplest of geometric objects, a point. However, $\mathbf{Spec}(\Z)$ is far from being a point, its underlying topological space is $0$ together with all prime numbers with the Zariski toplogy on them being the cofinite topology.

In the mid 1980-ties this led to the daydream that there might be an even more basic geometric object, the {\em absolute point} $\bullet$, and a map $\mathbf{Spec}(\Z) \rOnto \bullet$ such that $\mathbf{Spec}(\Z)$ might be viewed as a curve over $\bullet$. This absolute point would then be the geometric object corresponding to the {\em fake field} $\F_1$ with one element, and the hope was to mimic Weil's proof of the Riemann hypothesis for zeta-functions of curves over $\F_p$ to get at the Riemann hypothesis for prime numbers.

In these early days the mantra was that $\F_1$ forgets about the additive nature and only retains the multiplicative side of rings. That is, one would assume $\F_1$-algebras to be commutative monoids (or abelian groups), and the base-change map
\[
\Z \otimes_{\F_1} -~:~\mathbf{alg}_{\F_1} \rTo \mathbf{comm} \qquad M \mapsto \Z[M] \]
would assign to a monoid $M$ its integral monoid (or group) ring $\Z[M]$. In particular, the cyclic group of all $n$-th roots of unity $\pmb{\mu}_n$ would be interpreted as the extended fake field $\F_{1^n}$, with the algebraic closure of $\F_1$ then corresponding to $\pmb{\mu}_{\infty}$, the group of all roots of unity.

Let us try to understand which $\C$-algebraic varieties are defined over $\F_1$. A complex affine variety $X$ is defined over $\F_1$ if its coordinate ring $\Oscr(X)=\C[M]$ is a monoid algebra, so the free abelian monoid $\N^n$ gives us $\mathbb{A}^n_{\C}$ and the free abelian group $\Z^n$ gives us $T^n = \mathbb{G}^n_m$, the $n$-dimensional complex torus. Taking these as our building blocks one can then argue that the complex varieties defined over $\F_1$ are {\em torified varieties} as in \cite[\S 5.1]{LieberMM}. That is, $X = \sqcup_i T_i$ can be decomposed into tori $T_i$ and morphisms $X = \sqcup_i T_i \rTo Y = \sqcup_j T'_j$ are such that their restriction to every torus component is an algebraic morphism $T_i \rTo T'_{j(i)}$.

In $\mathbf{Mot}_{\C}$, the ring of motives of complex varieties, we have the {\em Lefschetz motive} $\mathbb{L}=[\mathbb{A}^1_{\C}]$ with $\mathbb{L}^n = [ \mathbb{A}^n_{\C}]$ and $[ \mathbb{G}_m^n ] = (\mathbb{L}-1)^n$. As a result, one can define as in \cite{LieberMM} the ring of motives of $\F_1$-varieties $\mathbf{Mot}_{\F_1}$ to be the subring $\Z[\mathbb{L}]$ of $\mathbf{Mot}_{\C}$ as any $\F_1$-variety $X$ decomposes over $\C$ into a disjoint union of tori, when its motive in $\mathbf{M}_{\C}$ is
\[
[ X ] = \sum_{i=0}^d a_i [\mathbb{G}_m^i ] = \sum_{i=0}^d a_i (\mathbb{L}-1)^i \]
where $a_i$ is the number of $i$-dimensional torus components. Lieber, Manin and Marcolli then define the counting measures
\[
\mu_{\F_{1^n}}~:~\mathbf{Mot}_{\F_1} = \Z[\mathbb{L}] \rTo \Z \qquad \mathbb{L} \mapsto 1+n \]
which count the number of points of an $\F_1$-variety $X = \sqcup_i T_i$ over $\F_{1^n}$, where $T_i$ is a $d(i)$-dimensional torus to be
\[
\# X(\F_{1^n}) = \sum_i n^{d(i)} \]
which is the number of points in the complex algebraic variety $X$ having all their coordinates in $\pmb{\mu}_n = \F_{1^n}$, see \cite[\S 5.5]{LieberMM}. Once again, one can then package these numbers in the {\em $\F_1$-zeta function}
\[
\zeta_{\F_1}(X) = exp(\sum_{n=1}^\infty \frac{\# X(\F_{1^n})}{n} t^n) \]
see \cite[\S 6.1]{LieberMM}.

\vskip 3mm

In this approach to $\F_1$-geometry one argues that any appearance of roots of unity from a geometric situation is an indication that there should be a corresponding variety over $\F_1$, as above or in any other implementation, encoding the essence of the situation. An illustrative example of this is due to Yu. Manin in \cite[\S 0.2]{Manin} and elaborated in \cite[\S 2]{MM}.

Consider a couple $(M,f)$ where $M$ is a compact manifold and $f$ is a {\em Morse-Smale diffeomorphism}, that is, $f$ is structurally stable and has a finite number of non-wandering points. Then, $f_*$ acts on the homology groups $H_k(M,\Z)$ (which are free $\Z$-modules) as a matrix $M(f_*)$. It is known that in this situation the eigenvalues of $M(f_*)$ are roots of unity. Manin argues that this action is similar to the action of the Frobenius on \'etale cohomology groups, in which case the eigenvalues are Weil numbers. That is, one might view roots of unity as Weil numbers in characteristic one. In \cite{Manin} he then asks to develop a version of $\F_1$-geometry, allowing for an object corresponding to the actions of $f_*$ on the $H_k(M,\Z)$. In \cite[\S 2.4.1]{MM} Manin and Marcolli propose such an object in the $\F_1$-geometry dreamed up by Jim Borger.

\section{Borger's idea}

So far, we considered base extension of varieties from $\F_1$ to $\Z$ or $\C$ and have seen that the integral or complex varieties defined over $\F_1$ are not especially interesting from a geometric viewpoint. We also didn't address the issue of viewing integral schemes, such as $\mathbf{Spec}(\Z)$, as geometric objects over $\F_1$, which is the problem of {\em forgetting the base} rather than of extending it.

Let us first consider the base extending/forgetting issue in the case of the finite field extension $\F_p \rInto \F_{p^n}$. Here, the {\em base extension} functor $v_* = - \times \F_{p^n}$
\[
\xymatrix{\mathbf{Var}_{\F_{p^n}} \ar@/^6ex/[dd]^{v^*} \ar@/_6ex/[dd]_{v_!} \\ \\ \mathbf{Var}_{\F_p} \ar[uu]^{v_*}} \]
has a right adjoint, which is{\em Weil descent}, and a left adjoint which is {\em forgetting the base}.

If we do have a suitable base extension functor $v_* : \mathbf{Var}_{\F_1} \rTo \mathbf{Var}_{\Z}$ we can make sense of the $\F_1$-variety corresponding to any integral scheme if this functor has a left adjoint. 

In the approach above, base extension $\mathbf{alg}_{\F_1} \rTo \mathbf{comm}$ was given by assigning to a commutative monoid $M$ its integral monoid ring $\Z[M]$. Alternatively, we can consider the sub-category $\mathbf{comm}^+_{mon}$ of $\mathbf{comm}$ consisting of all integral monoid rings and ring-morphisms coming from monoid maps. In this interpretation, the base extension functor is just the {\em forgetful functor} $F : \mathbf{comm}^+_{mon} \rTo \mathbf{comm}$. Further, as the direction of arrows reverses in going from rings to schemes, the ring-theoretical counter-part of a base forgetting functor should be a {\em right adjoint} to the forgetful functor. Unfortunately, whereas forgetful functors usually have a left adjoint (a universal construction) they seldom have a right adjoint (as is the case here).

Borger's idea, see \cite{Borger}, is to define the category of $\F_1$-algebras $\mathbf{alg}_{\F_1}$ to be a suitable sub-category $\mathbf{comm}_X^+$ of commutative rings with extra structure $X$ having the property that the forgetful functor $F$ not only has a left adjoint $U$ but also a right adjoint $G$
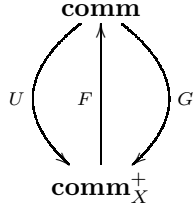
\[
\xymatrix{\mathbf{comm} \ar@/^6ex/[dd]^{G} \ar@/_6ex/[dd]_{U} \\ \\ \mathbf{comm}_X^+ \ar[uu]^{F}} \]
In such a situation, the coordinate ring of $\mathbf{Spec}(\Z)/\F_1$ is then $G(\Z) \in \mathbf{comm}_X^+=\mathbf{alg}_{\F_1}$ and the coordinate ring of the integral scheme corresponding to the {\em arithmatic plane} $\mathbf{Spec}(\Z) \times~_{\F_1} \mathbf{Spec}(\Z)$ is then $F\circ G(\Z)$, that is, stripping the extra structure $X$ from $G(\Z)$.

\vskip 3mm

In \cite{Borger}, Borger proposes to define $\mathbf{alg}_{\F_1}$ to be the sub-category $\mathbf{comm}_{\lambda}^+$ of commutative {\em $\lambda$-rings}. The motivation being that in case $A$ is a torsion-free $\Z$-ring, then $A$ is a $\lambda$-ring if and only if there is a commuting family $\{ \Psi_n~:~n \in \N \}$ of endomorphisms of $A$ with $\Psi_m \circ \Psi_n = \Psi_{mn}$ and such that for each prime number $p$ the endomorphism $\Psi_p$ is a lift of the Frobenius morphism on $A/pA$. So, the extra $\lambda$-ring structure can be interpreted as the {\em absolute Frobenius}. Further, the forgetful functor $F : \mathbf{comm}_{\lambda}^+ \rTo \mathbf{comm}$ does have a right adjoint, the functor $\mathbb{W}$ of {\em big Witt vectors}, see for example \cite{Hazewinkel}. In this proposal, the coordinate ring of $\mathbf{Spec}(\Z)/\F_1$ is
\[
\mathbb{W}(\Z) = 1 + t \Z [[t]] = \{ 1 + a_1 t + a_2 t^2 + \hdots~|~a_i \in \Z \} \subset \Z [[t]] \]
of all integral power series with constant term $1$ on which we define a new addition $\oplus$ which is ordinary multiplication of power series, and with a new multiplication $\otimes$ functorially induced by its action on geometric series
\[
\frac{1}{1-at} \otimes \frac{1}{1-bt} = \frac{1}{1-ab t} \]
and where the $\lambda$-ring structure is given by the endomorphisms $\Psi_n$ which are defined via their action on geometric series $\Psi_n(\tfrac{1}{1-at}) = \tfrac{1}{1-a^nt}$. A lot of subtle Galois theory is hidden in the phrase "functorially induced by the action on geometric series" in the definition of product and $\lambda$-ring structure.

\vskip 3mm

In \cite{LB} I proposed an alternative definition of $\mathbf{alg}_{\F_1}$, somewhat closer to the early days approach, but missing the absolute Frobenius feature. Here we take the sub-category $\mathbf{comm}_{bi}^+$ of torsion-free $\Z$-rings $A$ which are also birings, that is, have a comultiplication $\Delta : A \rTo A \otimes_{\Z} A$ and comultiplication $\epsilon : A \rTo \Z$. Note that integral monoid rings $\Z[M]$ are birings with $\Delta(m)=m \otimes m$ and $\epsilon(m)=1$ for all $m \in M$. In \cite{LB} it is proved that the forgetful functor $F : \mathbf{comm}_{bi}^+ \rTo \mathbf{comm}$ has a right adjoint $C$ given by taking the cocommutative free co-ring. In general, this is a horrible object, but fortunately in the case of interest to us, that is the coordinate ring of $\mathbf{Spec}(\Z)/\F_1$, we have that $C(\Z) = \mathbb{H}(\Z)$ the {\em Hadamard biring of integral linear recursive sequences}. That is,
\[
\mathbb{H}(\Z) = \{ (a_1,a_2,\hdots) \in \Z^{\infty}~|~\exists k,\exists B_1,..,B_k \in \Z~:~\forall n : a_n = a_{n-1}B_1+\hdots+a_{n-k}B_k \} \]
Addition and multiplication in $\mathbb{H}(\Z)$ comes from the component-wise operation in $\Z^{\infty}$. In contrast, the comultiplication is in general harder to describe and encodes subtle Galois information. If $\mathcal{M}$ is the multiplicative monoid of monic polynomials in $\Z[x]$, then the coring structure of $\mathbb{H}(\Z)$ is dual to the ring structures on quotients $\Z[x]/(F(x))$ for $F \in \mathcal{M}$ as
\[
\mathbb{H}(\Z) = \underset{\underset{F | G}{\rTo}}{lim}~(\frac{\Z[x]}{(F(x)})^* \]
and in particular $\epsilon(a_1,a_2,\hdots ) = a_1$. The sequence $d=(0,1,2,\hdots )$ is a primitive element in $\mathbb{H}(\Z)$, that is, $\Delta(d) = d \otimes 1 + 1 \otimes d$ and $\epsilon(d)=0$. One can use the structural result of commutative and co-commutative Hopf algebras over algebraically closed fields to deduce that 
\[
\mathbb{H}(\overline{\Q}) = (\overline{\Q}[\overline{\Q}^*_{\times}] \otimes \overline{\Q}[d])  \oplus K \]
with $\overline{\Q}[\overline{\Q}^*_{\times}]$ the group-algebra of the multiplicative group $\overline{\Q}^*$, $\overline{\Q}[d]$ is the enveloping algebra of the Lie-algebra $\overline{\Q} d$ and $K$ is the bialgebra-ideal of linear recursive sequences which are $0$ almost everywhere, see \cite{LarsonTaft}.

\section{Containers}

$\mathbb{W}(\Z)$ and $\mathbb{H}(\Z)$ are closer related than one might expect at first sight. In fact, we have a commuting square
\[
\xymatrix{\mathbb{W}_0(\Z) \ar[rr]^{L_{\Z}} \ar[d]_{Tr} & & \mathbb{W}(\Z) \ar[d]^{\mathghost} \\
\mathbb{H}(\Z) \ar[rr]_{i} & & \Z^{\infty}} \]
Here, $\mathbb{W}_0(\Z)$ is Almkvist's ring, see \cite{Almkvist}, constructed from pairs $(E,f)$ consisting of a projective (i.e. free) $\Z$-module with an endomorphism $f$. Such pairs are added resp. multiplied using direct sums resp. tensor products, so the zero-pair is $(0,0)$ and the one-pair is $(\Z,1)$. Almkvist's ring $\mathbb{W}_0(\Z)$ is then the quotient ring obtained by dividing out the ideal consisting of all pairs $(E,0)$. In fact, $\mathbb{W}_0$ is a functor on $\mathbf{comm}$ and one can show that $\mathbb{W}_0(\overline{\Q}) = \Z[\overline{\Q}^*_{\times}]$, so again a fair amount of Galois theory goes into the structure of $\mathbb{W}_0(\Z)$.

If $M_f$ is the integral matrix describing the endomorphism of a pair $(E,f)$, then we have a ringmorphism
\[
L_{\Z}~:~\mathbb{W}_0(\Z) \rTo \mathbb{W}(\Z) \qquad [E,f] \mapsto \frac{1}{det(1-tM_f)} \]
and in \cite[Thm. 6.4]{Almkvist} it is shown that the image of $L_{\Z}$ are precisely all {\em rational} integral formal power series. By taking the trace of the characteristic polynomial, we also have a map
\[
Tr~:~\mathbb{W}_0(\Z) \rTo \mathbb{H}(\Z) \qquad [E,f] \mapsto (Tr(M_f),Tr(M_f^2),\hdots ) \]
If we identify $t \Z[[t]]$ with $\Z^{\infty}$ via $\sum_{i \geq 1} a_i t^i \leftrightarrow (a_1,a_2,\hdots)$, the logarithmic derivative defines a ringmorphism (actually an isomorphism), called the {\em ghostmap}
\[
\mathghost~:~\mathbb{W}(\Z) \rTo \Z^{\infty} \qquad f(t) \mapsto t \frac{d}{dt} log(f(t)) \]
and the natural inclusion map $i$ turns this into a commuting square.

\vskip 3mm

We will argue that possible coordinate rings of $\mathbf{Spec}(\Z)/\F_1$, such as $\mathbb{W}(\Z)$ and $\mathbb{H}(\Z)$, should be viewed as {\em containers} for zeta-functions and other motivic information. So, let us return to geometric counting problems and the role these rings might play in them.

Often, one can equip a ring of motives $\mathbf{Mot}$ with a {\em pre $\lambda$-structure}, that is, a collection of maps $\lambda^n~:~\mathbf{Mot} \rTo \mathbf{Mot}$ for $n \in \N$ satisfying
\[
\lambda^0(X)=1,~ \quad \lambda^1(X)=X,~ \quad \text{and} \quad \lambda^n(X+Y) = \sum_{i+j=n} \lambda^i(X) \lambda^j(Y) \]
For example, the ring of motives of $\C$-varieties comes with a pre $\lambda$-structure given by $\lambda^n(X) = S^n(X)$, the $n$-th symmetric power of $X$. The pre $\lambda$-structure defines {\em Adams operations} $\{ \Psi^i \}$ for $i \geq 1$ via
\[
\mathghost(\sum_{i=0}^{\infty} \lambda^n(X) t^n) = (\Psi^1(X),\Psi^2(X),\hdots ) \]
where $\mathghost : \mathbb{W}(\mathbf{Mot}) \rTo \mathbf{Mot}^{\infty}$ is the ghostmap, that is the logarithmic derivative. It follows from the requirements on pre $\lambda$-structures that the Adams operations are additive. The ring of motives $\mathbf{Mot}$ is a $\lambda$-ring if and only if the Adams operations are also multiplicative.

Given a pre $\lambda$-structure $\{ \lambda^n \}$,one can associate to a counting measure $\mu : \mathbf{Mot} \rTo \Z$, {\em Kapranov's zeta function} corresponding to $\mu$
\[
\zeta_{\mu}(X) = \sum_{n=0}^{\infty} \mu(\lambda^n(X)) t^n \in 1 + t \Z[[t]] = \mathbb{W}(\Z) \]
One says that the counting measure $\mu$ is {\em exponentiable} if $\zeta_{\mu} : \mathbf{Mot} \rTo \mathbb{W}(\Z)$ is a ringmorphism, and it is said to be {\em rational} if there is a ringmorphism
$c_{\mu} : \mathbf{Mot} \rTo \mathbb{W}_0(\Z)$ such that $\zeta_{\mu} = L_{\Z} \circ c_{\mu}$.

For example, for the counting measure $\mu_{\F_p} : \mathbf{Mot}_{\F_p} \rTo \Z$ where $\mu_{\F_p}(X) = \#(X(\F_p))$ we have that Kapranov's zeta-function $\zeta_{\mu_{\F_p}}$ coincides with the Weil zeta-function $\zeta_W$ which is known to be rational. One can apply right-adjointness of the Witt functor $\mathbb{W}$ to show that the pre $\lambda$-structure given by $\lambda^n(X)=S^n(X)$ on $\mathbf{Mot}_{\C}$ does not equip $\mathbf{Mot}_{\C}$ with the structure of $\lambda$-ring, as this would imply that every {\em motivic measure}, that is a ringmorphism $\mu : \mathbf{Mot}_{\C} \rTo R$ with values in a commutative ring $R$, would be exponentiable, which is known to be not the case.

\vskip 3mm

Even if $\mathbf{Mot}$ is not a $\lambda$-ring, the pre $\lambda$-ring might induce a $\lambda$-ring structure on certain subrings of it. For example, it is known that the pre $\lambda$-structure $\lambda^n(X)=S^n(X)$ on $\mathbf{Mot}_{\C}$ makes the subring $\Z[\mathbb{L}]$ into a $\lambda$-ring such that $\Psi^n(\mathbb{L}) = \mathbb{L}^n$.
Thus, in the proposal where $\F_1$-varieties are torified varieties, and hence that $\mathbf{Mot}_{\F_1} = \Z[\mathbb{L}]$, it follows from right-adjointness, that is, the natural one-to-one correspondence
\[
\mathbf{comm}(\mathbf{Mot}_{\F_1},\Z) \leftrightarrow \mathbf{comm}^+_{\lambda}(\mathbf{Mot}_{\F_1},\mathbb{W}(\Z)) \qquad \mu \leftrightarrow \zeta_{\mu} \]
that every counting measure $\mu$ on $\mathbf{Mot}_{\F_1}$ is exponentiable, and even rational as there is a factorisation over $\mathbb{W}_0(\Z)$ defined by $\mathbb{L} \mapsto [ \Z,\mu(\mathbb{L}) ]$. Note that $\mathghost(\zeta_{\mu}(\mathbb{L}) = (\mu(\mathbb{L}),\mu(\mathbb{L})^2,\hdots )$.

However, the $\F_1$-zeta function of Lieber, Manin and Marcolli is {\em not} coming from a counting measure on $\mathbf{Mot}_{\F_1}$ in this way. For example, if $T^i$ denotes the $i$-dimensional torus $\mathbb{G}_m^i$ then
\[
\mathghost(\zeta_{\F_1}(T^i)) = t \frac{d}{dt} (\sum_{n=1}^{\infty} \frac{\# T^i(\F_{1^n}}{n}t^n) =  t \frac{d}{dt} (\sum_{n=1}^{\infty} n^{i-1} t^n) = (1,2^i,3^i,\hdots ) = (1,2,3,\hdots )^i \]
and observe that $(1,2,3,\hdots ) = 1 + d$ in $\Z^{\infty}$ (or in $\mathbb{H}(\Z)$). In particular, the $\F_1$-zeta function is not rational.

\vskip 3mm

What we did for counting measures with $\lambda$-rings and $\mathbb{W}(\Z)$ we can mimic for bi-rings and $\mathbb{H}(\Z)$. We say that a counting measure $\mu : \mathbf{Mot} \rTo \Z$ is {\em recursive} if it determines a ringmorphism $c_{\mu} : \mathbf{Mot} \rTo \mathbb{H}(\Z)$. If we are in a situation such that $\mathbf{Mot}$ can be given a bi-ring structure, it follows from right-adjointness that every counting measure on $\mathbf{Mot}$ is recursive. Once again, if we view the ring of motives of $\F_1$-varieties (in the torified interpretation) as the subring $\Z[\mathbb{L}]$ of $\mathbf{Mot}_{\C}$, then we can equip it with a bi-ring structure by demanding that $\mathbb{D}=\mathbb{L}-2$ is a primitive element, that is
\[
\Delta(\mathbb{D}) = \mathbb{D} \otimes 1 + 1 \otimes \mathbb{D} \quad \text{and} \quad \epsilon(\mathbb{D})=0 \]
By right-adjointness we have a natural one-to-one correspondence
\[
\mathbf{comm}(\mathbf{Mot}_{\F_1},\Z) \leftrightarrow \mathbf{comm}^+_{bi}(\mathbf{Mot}_{\F_1},\mathbb{H}(\Z)) \qquad\mu \leftrightarrow c_{\mu} \]
and hence every counting measure on $\mathbf{Mot}_{\F_1}$ is recursive. In particular, the counting measure $\mu_2$ defined by $\mu_2(\mathbb{L})=2$ is such that the corresponding bi-ring morphism $c_{\mu_2}$ satisfies
\[
\mathghost \circ \zeta_{\F_1} = c_{\mu_2} \]
as $c_{\mu_2}$ maps $\mathbb{D}$ to the primitive element $d \in \mathbb{H}(\Z)$.

\vskip 3mm

The commuting diagram above relating $\mathbb{W}_0(\Z)$ with $\mathbb{W}(\Z)$ and $\mathbb{H}(\Z)$ can also be used to give an answer to Manin's question of assigning objects in a variant of $\F_1$-geometry to a coupl;e $(M,f)$ where  $M$ is a compact manifold equipped with a Morse-Smale diffeomorphism $f$, see \cite[\S 0.2]{Manin}. We have seen that in such a situation the homology groups $H_k(M,\Z)$ with $0 \leq k \leq dim(M)$ are free $\Z$-modules on which $f_*$ acts as pre-multiplication with a matrix $M_k(f_*)$ having all its eigenvalues roots of unity.

To such a Morse-Smale couple $(M,f)$ we can therefore associate the family of elements $\{ [ H_k(M,\Z),M_k(f_*)]~:~0 \leq k \leq dim(M) \}$ in $\mathbb{W}_0(\Z)$. In the $\lambda$-ring variant of $\F_1$-geometry we can then associate the $\lambda$-subring $\mathbb{W}(M,f)$ of $\mathbb{W}(\Z)$ generated by the elements
\[
\{ L_{\Z}([ H_k(M,\Z),M_k(f_*)])~:~0 \leq k \leq dim(M) \} \]
and this is very similar to the approach by Manin and Marcolli in \cite[\S 2.4.1]{MM}. In the bi-ring approach to $\F_1$-geometry it is natural to associate to $(M,f)$ the sub bi-ring $\mathbb{H}(M,f)$ of the Hadamard ring $\mathbb{H}(\Z)$ generated by the elements
\[
\{ Tr([ H_k(M,\Z),M_k(f_*)])~:~0 \leq k \leq dim(M) \} \]
the hope being that these objects encode the essence of the dynamical system determined by the  Morse-Smale diffeomorphism.

\end{document}